\documentclass[12pt,a4paper,oneside,reqno]{article}

\usepackage[utf8]{inputenc}
\usepackage[english]{babel}
\usepackage[left=1in,right=1in,top=1in,bottom=1in]{geometry}
\usepackage{amsfonts,amssymb,amsmath,amsthm}
\usepackage[symbol]{footmisc}
\usepackage{longtable,array}
\newcolumntype{L}[1]{>{\raggedright\arraybackslash$}p{#1}<{$}}
\usepackage{enumerate}
\usepackage{cite}
\usepackage{xcolor}

\usepackage{hyperref}
\hypersetup{hidelinks}
\usepackage{fancyhdr}
\theoremstyle{plain}
\newtheorem{theorem}{Theorem}

\newtheorem{corollary}[theorem]{Corollary}
\theoremstyle{definition}

\makeatletter
\@removefromreset{theorem}{section}\@removefromreset{equation}{section}

\makeatother
\usepackage{adjustbox}
\usepackage{tikz}\usetikzlibrary{arrows}

\DeclareMathOperator{\Id}{Id}
\DeclareMathOperator{\Aut}{Aut}

\IfFileExists{fouriernc.sty}{\usepackage{fouriernc}}{}

\newcommand{\doi}[1]{\href{https://doi.org/#1}{doi:\detokenize{#1}}}
\newcommand{\arxiv}[1]{\href{https://arxiv.org/abs/#1}{arXiv:#1}}

\theoremstyle{plain}
\newtheorem{Th}[theorem]{Theorem}
\newtheorem{Lem}[theorem]{Lemma}
\theoremstyle{definition}
\newtheorem{Remark}[theorem]{Remark}

\newcommand{\FF}{{\mathbb{F}}}

\newcommand{\OO}{\mathcal{O}}
\newcommand{\SSS}{\mathcal{S}}
\newcommand{\TT}{\mathcal{T}}

\begin{document}

\noindent\fbox{\parbox{\dimexpr\textwidth-2\fboxsep-2\fboxrule\relax}{\small
\textbf{Corrected version (verified by Claude).} This is a corrected version of the paper I.~Kaygorodov and Yu.~Volkov, \emph{Conservative algebras of $2$-dimensional algebras, II}, Communications in Algebra, 45 (2017), no.~8, 3413--3421, \doi{10.1080/00927872.2016.1236935}.}}

\bigskip
\bigskip
\begin{center}
{\Large Conservative algebras of $2$-dimensional algebras, II}

\bigskip
\bigskip
{\bf Ivan Kaygorodov and Yury Volkov}
\end{center}

\bigskip

\noindent{\small\textbf{Abstract.} In 1990 Kantor defined the conservative algebra $W(n)$ of all algebras (i.e. bilinear maps) on the $n$-dimensional vector space. If $n>1$, then the algebra $W(n)$ does not belong to any well-known class of algebras (such as associative, Lie, Jordan, or Leibniz algebras). We describe automorphisms, one-sided ideals, and idempotents of $W(2)$. Similar problems are also solved for the algebra $W_2$ of all commutative algebras on the $2$-dimensional vector space and for the algebra $S_2$ of all commutative algebras with trace zero multiplication on the $2$-dimensional vector space.}

\bigskip

\bigskip
\bigskip

\bigskip
\bigskip

\tableofcontents

\bigskip

\newpage
\section*{Introduction}\addcontentsline{toc}{section}{Introduction}
Throughout this paper, $\FF$ is a fixed field of characteristic zero.
Algebras under consideration are not necessarily
unital or associative.
A multiplication on a vector space $W$ is a bilinear mapping $W\times W \to W$. We denote by $(W,P)$ the algebra with underlying space $W$ and multiplication $P$. Given a vector space $W$, a linear mapping $A:W\rightarrow W$, and a bilinear mapping $B:W\times W\to W$, we can define a multiplication $[A,B]:W\times W\to W$ by the formula
$$[A,B](x,y)=A(B(x,y))-B(A(x),y)-B(x,A(y))$$
for $x,y\in W$. For an algebra $A$ with a multiplication $P$ and $x\in A$ we denote by $L_x^P$ the operator of left multiplication by $x$. If the multiplication $P$ is fixed, we write $L_x$ instead of $L_x^P$.

In 1990 Kantor~\cite{Kantor90} defined the multiplication $\cdot$ on the set of all algebras (i.e. all multiplications) on the $n$-dimensional vector space $V_n$ as follows:
$$A\cdot B = [L_{e}^A,B],$$
where $A$ and $B$ are multiplications and $e\in V_n$ is some fixed vector.
Let $W(n)$  denote the algebra of all algebra structures on $V_n$ with multiplication defined above.
If $n>1$, then the algebra $W(n)$ does not belong to any well-known class of algebras (such as associative, Lie, Jordan, or Leibniz algebras). The algebra $W(n)$ turns out to be a conservative algebra (see below).

In 1972 Kantor~\cite{Kantor72} introduced conservative algebras as a generalization of Jordan algebras. Namely, an algebra $A=(W,P)$ is called a {\it conservative algebra} if there is a new multiplication $F:W\times W\rightarrow W$ such that
\begin{eqnarray}\label{tojd_oper}
[L_b^P, [L_a^P,P]]=-[L_{F(a,b)}^P,P]
\end{eqnarray}
for all $a,b\in W$. In other words, the following identity holds for all $a,b,x,y\in W$:
\begin{multline}\label{tojdestvo_glavnoe}
b(a(xy)-(ax)y-x(ay))-a((bx)y)+(a(bx))y+(bx)(ay)\\
-a(x(by))+(ax)(by)+x(a(by))
=-F(a,b)(xy)+(F(a,b)x)y+x(F(a,b)y).
\end{multline}
The algebra $(W,F)$ is called an algebra {\it associated} to $A$.

Let us recall some well-known results about conservative algebras. In~\cite{Kantor72} Kantor  classified all simple conservative algebras and triple systems of second order and defined the class of  terminal algebras as algebras satisfying a certain identity. He proved that every terminal algebra is a conservative algebra and classified all simple finite-dimensional terminal algebras with a left quasi-unit over an algebraically closed field of characteristic zero~\cite{Kantor89term}. Terminal trilinear operations were studied in~\cite{Kantor89tril}, and some questions concerning classification of simple conservative algebras were considered in~\cite{Kantor89trudy}. After that, Cantarini and Kac classified simple finite-dimensional (and linearly compact) super-commutative and super-anticommutative conservative superalgebras and some generalization of these algebras (also known as ``rigid''{} superalgebras) over an algebraically closed field of characteristic zero \cite{kac_can_10}.

The algebra $W(n)$ plays a similar role in the theory of conservative algebras as the Lie algebra of all $n\times n$ matrices $gl_n$ plays in the theory of Lie algebras. Namely, in~\cite{Kantor88,Kantor90}  Kantor considered the category $\mathcal{S}_n$ whose objects are conservative algebras of non-Jacobi dimension $n$. It was proven that the algebra $W(n)$ is the universal attracting object in this category, i.e., for every $M\in\mathcal{S}_n$ there exists a canonical homomorphism from $M$ into the algebra $W(n)$. In particular, all Jordan algebras of dimension $n$ with unity are contained in the algebra
$W(n)$. The same statement also holds for all noncommutative Jordan algebras of dimension $n$ with unity.
Some properties of the product in the algebra $W(n)$ were studied in \cite{kay}.
The study of low dimensional conservative algebras was started in \cite{kaylopo},
where derivations and subalgebras of codimension 1 of the algebra $W(2)$ and of its principal subalgebras $W_2$ and $S_2$ were described.
Note that $W_2$ and $S_2$ are simple terminal algebras with a left quasi-unit from the classification of Kantor~\cite{Kantor89term}.

\section{Conservative algebra $W(2)$ and its subalgebras}

A multiplication on a $2$-dimensional vector space is defined by a $2\times 2\times 2$ array of structure constants. The classification of such multiplications was given in many papers (see, for example,~\cite{GR11}).
Let us consider the space $W(2)$ of all multiplications on the 2-dimensional space $V_2$ with a basis $v_1,v_2$.  The definition of the multiplication $\cdot$  on the algebra $W(2)$ can be found in  Section~1. Namely, we fix the vector $v_1 \in V_2$ and define
$$(A \cdot B)(x,y)=A(v_1,B(x,y))-B(A(v_1,x),y)-B(x,A(v_1,y))$$
for $x,y \in V_2$ and $A,B \in W(2)$. The algebra $W(2)$ is conservative \cite{Kantor90}.

Let us consider the multiplications $\alpha_{i,j}^k$ ($i,j,k=1,2$) on $V_2$ defined by the formula $\alpha_{i,j}^k(v_t,v_l)=\delta_{it}\delta_{jl} v_k$ for all $t,l$. It is easy to see that $\{ \alpha_{i,j}^k | i,j,k=1,2 \}$ is a basis of the algebra $W(2)$. The multiplication table of $W(2)$ in this basis is given in \cite{kaylopo}. In this work we use another basis for the algebra $W(2)$. Let us introduce the notation
$$
\begin{aligned}
e_1&=\alpha_{11}^1-\alpha_{12}^2-\alpha_{21}^2,
e_2=\alpha_{11}^2, e_3=\alpha_{22}^2-\alpha_{12}^1-\alpha_{21}^1,
e_4=\alpha_{22}^1,\\ e_5&=2\alpha_{11}^1+\alpha_{12}^2+\alpha_{21}^2,
e_6=2\alpha_{22}^2+\alpha_{12}^1+\alpha_{21}^1,
e_7=\alpha_{12}^1-\alpha_{21}^1,
e_8=\alpha_{12}^2-\alpha_{21}^2.
\end{aligned}
$$
It is easy to see that the multiplication table of $W(2)$ in the basis $e_1,\dots,e_8$ is the following.

\begin{center}
\begin{tabular}{c|c|c|c|c|c|c|c|c|}
      & $e_1$ & $e_2$ & $e_3$ & $e_4$ & $e_5$ & $e_6$ & $e_7$ & $e_8$ \\ \hline
$e_1$ & $-e_1$ & $-3e_2$ & $e_3$ & $3e_4$ & $-e_5$ & $e_6$ & $e_7$ & $-e_8$ \\ \hline
$e_2$ & $3e_2$ & $0$ & $2e_1$ & $e_3$ & $0$ & $-e_5$ & $e_8$ & $0$ \\ \hline
$e_3$ & $-2e_3$ & $-e_1$ & $-3e_4$ & $0$ & $e_6$ & $0$ & $0$ & $-e_7$ \\ \hline
$e_4$ & $0$ & $0$ & $0$ & $0$ & $0$ & $0$ & $0$ & $0$ \\ \hline
$e_5$ & $-2e_1$ & $-3e_2$ & $-e_3$ & $0$ & $-2e_5$ & $-e_6$ & $-e_7$ & $-2e_8$ \\ \hline
$e_6$ & $2e_3$ & $e_1$ & $3e_4$ & $0$ & $-e_6$ & $0$ & $0$ & $e_7$ \\ \hline
$e_7$ & $2e_3$ & $e_1$ & $3e_4$ & $0$ & $-e_6$ & $0$ & $0$ & $e_7$ \\ \hline
$e_8$ & $0$ & $e_2$ & $-e_3$ & $-2e_4$ & $0$ & $-e_6$ & $-e_7$ & $0$ \\ \hline
\end{tabular}
\end{center}
\medskip
The subalgebra spanned by the elements $e_1, \ldots, e_6$ is the conservative (and, moreover,  terminal) algebra $W_2$ of commutative 2-dimensional algebras.
The subalgebra spanned by the elements $e_1, \ldots, e_4$ is the conservative (and, moreover,  terminal) algebra $S_2$ of all commutative 2-dimensional algebras with trace zero multiplication \cite{kaylopo}.
We denote by $p_k:W(2)\rightarrow\FF$ ($k=1,\dots,8$) the map which sends $\sum\limits_{k=1}^8a_ke_k$ to $a_k$.

\section{One-sided ideals}

Let $Ann_l(A)$ denote the left annihilator of the algebra $A$.
It is easy to see that  $Ann_l(W(2))$ is spanned by the elements
$e_4, e_5-2e_1-3e_8, e_3+e_6$ and $e_3+e_7,$
  $Ann_l(W_2)$ is spanned by the elements $e_4$ and $e_3+e_6,$ and
  $Ann_l(S_2)$ is spanned by the element $e_4.$

\begin{Lem}\label{righideal}
If $I$ is a non-trivial right ideal of $W(2)$,
then $I$ is a subspace of $Ann_l(W(2)).$
\end{Lem}

{\bf Proof.}
Note that if $e_1\in I$, then $I=W(2)$.
Let $x\in I.$ If $p_2(x) \neq 0$, then
$(xe_4)e_2 = -p_2(x) e_1 \in I.$
In this case $e_1\in I$ and $I=W(2)$.
Let us now consider the case where $p_2(x)=0$ for any $x\in I$.

We have $(x e_3)e_2=- (p_1-p_5-p_8)(x)e_1 \in I$.
So $I=W(2)$ or $p_1(x)=p_5(x)+p_8(x)$ for any $x\in I$.
We have
$$xe_2= (-p_3+p_6+p_7)(x)e_1+(-3p_1-3p_5+p_8)(x)e_2.$$
Since $xe_2\in I$, we have $p_8(x)=3(p_1+p_5)(x)$. Then $I=W(2)$ or $p_3(x)=(p_6+p_7)(x)$.
So, if $I\neq W(2)$, then any element $x\in I$ has the form
$$x=p_5(x)(e_5-2e_1-3e_8)+p_4(x)e_4+p_6(x)(e_3+e_6)+p_7(x)(e_3+e_7)\in Ann_l(W(2)).$$
The lemma is proved.

Similarly, one can prove the following two lemmas.

\begin{Lem}
If $I$ is a non-trivial right ideal of $W_2$,
then $I$ is a subspace of  $Ann_l(W_2).$
\end{Lem}

\begin{Lem}
If $I$ is a non-trivial right ideal of $S_2$, then $I$ is a subspace of  $Ann_l(S_2).$
\end{Lem}

Now we describe left ideals of the algebras $W(2),W_2$ and $S_2.$

\begin{Lem}\label{leftideal}
If $I$ is a non-trivial left  ideal of $W(2)$,
then $I$ is one of the following spaces:

1) $W_{\alpha,\beta}={}_{\FF}\langle e_1,e_2,e_3,e_4,  \alpha e_5 +\beta e_8, -\alpha e_6 + \beta e_7\rangle$, where
$\alpha, \beta \in \FF$, $(\alpha,\beta)\not=(0,0)$;

2) $I_1={}_{\FF}\langle e_1,e_2,e_3,e_4\rangle;$

3) $I_2={}_{\FF}\langle e_5,e_6,e_7,e_8\rangle;$

4) $w_{\alpha,\beta}={}_{\FF}\langle\alpha e_5 +\beta e_8, -\alpha e_6 + \beta e_7\rangle,$
where $\alpha, \beta \in \FF$, $(\alpha,\beta)\not=(0,0)$.

\end{Lem}

{\bf Proof.}
Let $x\in I.$ We have
$e_2(e_2(e_1x-x))= 4 p_4(x) e_1 \in I.$

{\bf 1.} Assume that there is $x\in I$ such that $p_4(x) \neq 0.$
Then it is obvious that $e_1,e_2,e_3,e_4 \in I.$
In this case either $I$ is the subspace generated by $e_1,e_2,e_3,e_4$ (in this case $I=I_1$) or
there exists a nonzero element
$w=\sum\limits_{i=5}^8 \alpha_i e_i  \in I.$
In the latter case we have
\begin{eqnarray*}
 w_1=w-e_1w &=& 2\alpha_5e_5+2\alpha_8e_8\in I,\\
 w_2=e_2(e_1w+w)&=&-2\alpha_6e_5+2\alpha_7e_8\in I,\\
 w_3=e_1w+w &=& 2\alpha_6e_6+2\alpha_7e_7 \in I,\\
 w_4=-e_3(w-e_1w)&=&-2\alpha_5e_6+2\alpha_8e_7 \in I.\\
 \end{eqnarray*}
Now, if the elements $w_1$ and $w_2$ are linearly independent, then $e_5,e_6,e_7,e_8 \in I$ and $I=W(2)$.
If the elements $w_1$ and $w_2$ are linearly dependent for every such $w$, then $I=W_{\alpha,\beta}$ for some $\alpha,\beta\in\FF$.

{\bf 2.} Assume that $p_4(x)=0$ for any $x\in I$.
If $p_i(x)=0$ for any $x\in I$ and $i=1,2,3,4,$
then by the previous argument, we have either $I=I_2$ or $I=w_{\alpha,\beta}$ for some $\alpha,\beta\in\FF$.

Let us consider $x\in I$. If $p_1(x)\neq 0$, then we have $$e_6(e_3(e_3(e_2(e_1x-x))))=-36 p_1(x) e_4 \in I$$
and we have a contradiction with $p_4(I)=0.$

Suppose that $p_1(I)=p_4(I)=0$.
If $p_2(x)\neq 0$ or $p_3(x)\neq0$ for some $x\in I$, then
we have $$e_3(e_1x+x)=2p_2(x)e_1 -6p_3(x) e_4 \in I$$
and we have a contradiction with $p_1(I)=p_4(I)=0$.
The lemma is proved.

Similarly, one can prove the following two lemmas.

\begin{Lem}\label{leftW2}
If $I$ is a non-trivial left  ideal of $W_2$,
then $I$ is one of the following spaces:

1) ${}_{\FF}\langle e_1,e_2,e_3,e_4\rangle;$

2) ${}_{\FF}\langle e_5,e_6\rangle.$

\end{Lem}

\begin{Lem}
There are no non-trivial left ideals in $S_2$.
\end{Lem}

\subsection{Corollaries}
Here we want to formulate some corollaries about a decomposition of the algebras under consideration in a sum of one-sided ideals and about their ternary derivations.

\begin{corollary} For any left ideal $Y_1$ of the algebra $W \in \{W(2), W_2,S_2\}$ there is some left ideal $Y_2$
of $W$ such that $W=Y_1+Y_2$ and $Y_1 \cap Y_2=0.$
\end{corollary}

The definition of a ternary derivation arose in the paper of Jim\'enez-Gestal and P\'erez-Izquierdo~\cite{perez03}.
A ternary derivation is a triple of linear mappings $(D,F,G)$ such that $D(xy)=F(x)y+xG(y).$
A ternary derivation is trivial if $D$, $F$, and $G$ are sums of a derivation and an element of a centroid of an algebra under consideration. Ternary derivations of Jordan algebras and superalgebras, and $n$-ary algebras were studied in \cite{kay13,kay14,shestakov14}.

\begin{corollary} The algebras $W(2), W_2$ and $S_2$ have non-trivial ternary derivations.
\end{corollary}

{\bf Proof.} Let $A$ be an algebra from the set $\{W(2),W_2,S_2\}.$
By the lemmas above $A$ has a nonzero left annihilator $Ann_l(A)$.
Then for any nonzero linear map $\phi$ such that $\phi(A) \subseteq Ann_l(A)$  the triple $(0, \phi, 0)$ is a non-trivial ternary derivation of $A$.
The corollary is proved.

\medskip

\section{Automorphisms}

In this section we discuss the automorphism groups of the algebras $W(2)$, $W_2$ and $S_2$. Firstly we prove a lemma which provides a series of automorphisms of the algebra $W(n)$ for any $n\ge 2$. Let us fix a space $V$ and a nonzero element $e\in V$. We denote by $W(V,e)$ the algebra whose elements are bilinear maps from $V\times V$ to $V$ and multiplication on which is defined by the formula
$$
(A\cdot B)(x,y)=A(e,B(x,y))-B(A(e,x),y)-B(x,A(e,y))
$$
for $A,B\in W(V,e)$, $x,y\in V$. Let us introduce the notation $GL(V,e)=\{f\in GL(V)\mid f(e)=e\}$. It is clear that $GL(V,e)$ is a subgroup of $GL(V)$. If $V$ is an $n$-dimensional vector space, then $GL(V,e)$ is isomorphic to $(n-1)$-dimensional {\it affine group} over $\FF$.

\begin{Lem}\label{auttoaut}
There is a monomorphism of groups $\Phi:GL(V,e)\rightarrow \Aut\big(W(V,e)\big)$ defined by the equality
$$
\big(\Phi(f)(A)\big)(x,y)=fA(f^{-1}(x),f^{-1}(y))
$$
for $f\in GL(V,e)$, $A\in W(V,e)$, $x,y\in V.$
\end{Lem}
{\bf Proof.}
Let us prove that $\Phi(f)\in \Aut\big(W(V,e)\big)$ for any $f\in GL(V,e)$. Actually, for $A,B\in W(V,e)$, $x,y\in V$ we have
\begin{multline*}
\big(\Phi(f)(A)\cdot \Phi(f)(B)\big)(x,y)=
fA\big(f^{-1}(e),f^{-1}fB(f^{-1}(x),f^{-1}(y))\big)\\
-fB\big(f^{-1}fA(f^{-1}(e),f^{-1}(x)),f^{-1}(y)\big)
-fB\big(f^{-1}(x),f^{-1}fA(f^{-1}(e),f^{-1}(y))\big)\\
=f\Big(A\big(e,B(f^{-1}(x),f^{-1}(y))\big)
-B\big(A(e,f^{-1}(x)),f^{-1}(y)\big)
-B\big(f^{-1}(x),A(e,f^{-1}(y))\big)\Big)\\
=\big(\Phi(f)(A\cdot B)\big)(x,y).
\end{multline*}
It is easy to check that $\Phi$ is a homomorphism of groups.
Let $f\in GL(V,e)$ be such that $f\not=\Id_V$. Then there is such $x\in V$ that $f^{-1}(x)\not=x$. It is clear that there is such $A\in W(V,e)$ that $A(e,x)=e$ and $A(e,f^{-1}(x))\not=e$. Then
$$\Phi(f)(A)(e,x)=fA(f^{-1}(e),f^{-1}(x))=fA(e,f^{-1}(x))\not=f(e)=A(e,x).$$
So $\Phi(f)\not=\Id_{W(V,e)}$ and we prove that $\Phi$ is a monomorphism.
The lemma is proved.

\medskip

Any element of the group $GL(V_2,v_1)$ can be written as a composition $\SSS_b\TT_a=\TT_{ab}\SSS_b$ ($a,b\in\FF$, $b\not=0$), where $\TT_a(v_2)=v_2+av_1$ and $\SSS_b(v_2)=\frac{v_2}{b}$.
We write simply $\TT_a$ and $\SSS_b$ instead of $\Phi(\TT_a)$ and $\Phi(\SSS_b)$. Direct calculations show that the action of $\TT_a$ and $\SSS_b$ on $W(V_2)$ is defined by the formulas
\begin{multline*}
\TT_a(e_1)=e_1+2ae_3+3a^2e_4,\:\:\TT_a(e_2)=e_2+ae_1+a^2e_3+a^3e_4,\:\:\TT_a(e_3)=e_3+3ae_4,\:\:\TT_a(e_4)=e_4,\\
\TT_a(e_5)=e_5-ae_6,\:\:\TT_a(e_6)=e_6,\:\:\TT_a(e_7)=e_7,\:\:\TT_a(e_8)=e_8+ae_7
\end{multline*}
and
\begin{multline*}
\SSS_b(e_1)=e_1,\:\:\SSS_b(e_2)=\frac{e_2}{b},\:\:
\SSS_b(e_3)=be_3,\:\:\SSS_b(e_4)=b^2e_4,\\
\SSS_b(e_5)=e_5,\:\:\SSS_b(e_6)=be_6,\:\:
\SSS_b(e_7)=be_7,\:\:\SSS_b(e_8)=e_8.
\end{multline*}
It follows from these formulas that $\TT_a$ and $\SSS_b$ induce automorphisms of $W_2$ and $\SSS_2$ which we denote by $\TT_a$ and $\SSS_b$ too.

\begin{Remark}\label{rem_exp}
It follows from the description of the derivations of $W(2)$ in~\cite{kaylopo} that $L_{e_7}$ and $L_{e_8}$ are derivations of the algebra $W(2)$. The map $L_{e_7}$ is nilpotent and it is easy to see that $\TT_a=\sum_{k=0}^{\infty}\frac{(aL_{e_7})^k}{k!}$.
On the other hand the map $L_{e_8}$ is not nilpotent. But if the element $e^t=\sum_{k=0}^{\infty}\frac{t^k}{k!}$ is defined and lies in $\FF^*$ for some $t\in\FF$, then we have an equality $\SSS_{e^{-t}}=\sum_{k=0}^{\infty}\frac{(tL_{e_8})^k}{k!}$. So if for any $b\in\FF^*$ there is $t\in\FF$ such that $b=e^t$, then the automorphism $\SSS_b$ can be defined using  the derivation $L_{e_8}$.
\end{Remark}

Now we are ready to describe the automorphism groups of the algebras $W(2)$, $W_2$ and $S_2$.

\begin{Th}
The homomorphism $\Phi:GL(V_2,v_1)\rightarrow \Aut\big(W(2)\big)$ is bijective. Moreover, the induced homomorphisms from $GL(V_2,v_1)$ to $\Aut(W_2)$ and $\Aut(S_2)$ are bijective too. In other words, $\Aut\big(W(2)\big)\cong \Aut(W_2)\cong \Aut(S_2)$ is isomorphic to the matrix group $\begin{pmatrix}1 & \FF\\ 0 & \FF^*\end{pmatrix}$ (the one-dimensional affine group over $\FF$).
\end{Th}
{\bf Proof.} It follows from the proof of Lemma \ref{auttoaut} that all three homomorphisms of groups considered in the theorem are monomorphic. So it is enough to prove that they are epimorphic. Firstly we prove that the automorphisms $\TT_a$ ($a\in\FF$) and $\SSS_b$ ($b\in\FF^*$) generate the group $\Aut(S_2)$.

Let us take some $f\in\Aut(S_2)$. Since $f$ sends any right ideal to a right ideal, we have $f(e_4)=te_4$ for some $t\in\FF^*$. Then we have $0=f(e_3e_4)=tf(e_3)e_4$. So $f(e_3)=be_3+ae_4$ for some $a,b\in\FF$. Let us consider the equality $-3f(e_4)=f(e_3)f(e_3)$. It can be rewritten in the form $-3te_4=-3b^2e_4$. So $t=b^2$. Let us consider $g=\SSS_{\frac{1}{b}}\TT_{-\frac{a}{3b}}f$. We know that $g(e_3)=e_3$ and $g(e_4)=e_4$ and it remains to prove that $g=\Id_{S_2}$. It follows from the equality $e_3=g(e_1e_3)=g(e_1)e_3$ that $g(e_1)=e_1+te_4$ for some $t\in\FF$. Then we have $$-e_1-te_4=g(-e_1)=g(e_1e_1)=\big(g(e_1)\big)^2=-e_1+3te_4.$$
So $t=0$ and we have $g(e_1)=e_1$. Using the equality $2e_1=g(e_2e_3)=g(e_2)e_3$ we obtain that $g(e_2)=e_2+te_4$ for some $t\in\FF$. Consideration of the equality $0=\big(g(e_2)\big)^2$ shows that $t=0$ and $g(e_2)=e_2$. So we prove that $\Aut(S_2)$ is generated by $\TT_a$ and $\SSS_b$ ($a,b\in\FF$, $b\not=0$).

Let us now prove that $\TT_a$ ($a\in\FF$) and $\SSS_b$ ($b\in\FF^*$) generate the group $\Aut(W_2)$. Since $S_2$ is the unique 4-dimensional left ideal of $W_2$ by Lemma \ref{leftW2}, every automorphism of $W_2$ induces an automorphism of $S_2$. Let us take some $f\in\Aut(W_2)$. Since $\Aut(S_2)$ is generated by automorphisms of the form $\TT_a$ and $\SSS_b$, there are such $a\in\FF$ and $b\in\FF^*$ that the automorphism $g=\SSS_b\TT_af$ is identical on $S_2$. It remains to prove that $g(e_5)=e_5$ and $g(e_6)=e_6$.
But $e_5$ and $e_6$ generate the unique 2-dimensional left ideal of $W_2$ and so $g(e_5)=t_1e_5+t_2e_6$, $g(e_6)=t_3e_5+t_4e_6$. Consideration of the equalities $-3e_2=g(e_5e_2)=g(e_5)e_2$ and $e_1=g(e_6e_2)=g(e_6)e_2$ shows that $t_1=t_4=1$ and $t_2=t_3=0$, i.e. $g=\Id_{W_2}$.

It remains to prove that $\TT_a$ ($a\in\FF$) and $\SSS_b$ ($b\in\FF^*$) generate the group $\Aut\big(W(2)\big)$. By Lemma \ref{leftideal} the algebra $W(2)$ has two 4-dimensional left ideals: $S_2$ and $I={}_{\FF}\langle e_5,e_6,e_7,e_8\rangle$. Note that $W(2)$ contains $2$-dimensional left ideals and they all lie in $I$. So any automorphism of the algebra $W(2)$ maps $S_2$ to $S_2$ and $I$ to $I$. Analogously to the case of the algebra $W_2$ it is enough to prove that any automorphism $g\in\Aut\big(W(2)\big)$ identical on $S_2$ is identical on $I$. It follows from the equalities $g(e_5)=-e_1g(e_5)$ and $g(e_8)=-e_1g(e_8)$ that the space ${}_{\FF}\langle e_5,e_8\rangle$ is invariant under the action of $g$. Using the equalities $-2e_1=g(e_5)e_1$ and $0=g(e_5)e_4$ we see that $g(e_5)=e_5$. Using the equalities $0=g(e_8)e_1$ and $-2e_4=g(e_8)e_4$ we see that $g(e_8)=e_8$. Finally, we have $g(e_6)=g(e_3e_5)=e_3e_5=e_6$ and $g(e_7)=-g(e_3e_8)=-e_3e_8=e_7$. So the theorem is proved.

\begin{corollary}\label{isoaut}
Let $A$ and $B$ be two multiplications on the space $V_2$. There is an isomorphism from the algebra $(V_2,A)$ to the algebra $(V_2,B)$ which sends $v_1$ to $v_1$ iff there is an automorphism of the algebra $W(2)$ which sends $(V_2,A)$ to $(V_2,B)$.
\end{corollary}

\section{Idempotents}

Let us now describe the idempotents of the algebra $W(2)$. For $e\in W(2)$ we denote by $\OO(e)=\{\SSS_b\TT_a(e)\mid a,b\in\FF,b\not=0\}$ the orbit of $e$ under the action of the automorphism group of $W(2)$. It is clear that $e$ is an idempotent iff all elements of $\OO(e)$ are idempotents.
We denote by $\bar\FF$ the set of representatives of cosets of the subgroup $(\FF^*)^2$ in the group $\FF^*$.

\begin{Th}\label{idemp}
The set of nonzero idempotents of the algebra $W(2)$ is the disjoint union of the following sets:
\begin{enumerate}
\item $\OO(e_8+e_2-e_1+c(3e_8+e_5-2e_1))$ ($c\in\FF$);
\item $\OO(-e_1+c(e_5-2e_1)+de_8)$ ($c,d\in\FF$);
\item $\OO(-e_1-2e_8+4e_3+e_6+3e_7+c(3e_8-e_5+2e_1)+de_4)$ ($c,d\in\FF$);
\item $\OO(-e_1-2e_8+c(3e_8-e_5+2e_1)+qe_4)$ ($c\in\FF$, $q\in\bar\FF$);
\end{enumerate}
\end{Th}
{\bf Proof.} Direct calculations show that all elements of the set described in the theorem are idempotent. Let $e\in W(2)$ be a nonzero idempotent. Let us prove that $e$ lies in the set described in the theorem. For $c,d\in \FF$, $q\in\bar\FF$ let us introduce the notation
\begin{multline*}
w_1(c)=e_8+e_2-e_1+c(3e_8+e_5-2e_1),
w_2(c,d)=-e_1+c(e_5-2e_1)+de_8,\\
w_3(c,d)=-e_1-2e_8+4e_3+e_6+3e_7+c(3e_8-e_5+2e_1)+de_4,\\
w_4(c,d)=-e_1-2e_8+c(3e_8-e_5+2e_1)+de_4.
\end{multline*}

{\bf 1.} Assume that $p_2(e)\not=0$. Let us take
$b=p_2(e)$ and $a=-\frac{p_1(e)+2p_5(e)+1}{p_2(e)}$.
It is easy to see that $p_2\SSS_b\TT_a(e)=1$ and $(p_1+2p_5)\SSS_b\TT_a(e)=-1$. So we may suppose that $p_2(e)=1$ and $(p_1+2p_5)(e)=-1$ initially.
Let us introduce the notation $c=p_5(e)$. Then we have $p_1(e)=-(2c+1)$. It follows from the equality $p_2(e^2)=p_2(e)$ that $p_8(e)=3c+1$. Rewriting the equalities $p_5(e^2)=p_5(e)$ and $p_8(e^2)=p_8(e)$ we obtain that $p_6(e)=p_7(e)=0$. Rewriting the equality $p_1(e^2)=p_1(e)$ we obtain $p_3(e)=-p_6(e)-p_7(e)=0$. Then we have $0=p_3(e^2)=p_2(e)p_4(e)=p_4(e)$. So we have
$$
e=\sum\limits_{k=1}^8p_k(e)e_k=-(2c+1)e_1+e_2+ce_5+(3c+1)e_8=w_1(c).
$$

Let us now prove that $(p_1+2p_5)(e)=-1$ if $e$ is a nonzero idempotent such that $p_2(e)=0$. Using the multiplication table we see that
$$p_1(e^2)=-\big(p_1(e)+2p_5(e)\big)p_1(e),\:\:
p_5(e^2)=-\big(p_1(e)+2p_5(e)\big)p_5(e).$$
Since $e^2=e$, we have $p_1(e)+2p_5(e)=-1$ in the case where $p_1(e)\not=0$ or $p_5(e)\not=0$. It is enough to prove that $e=0$ if $p_1(e)=p_2(e)=p_5(e)=0$. But in this case we obtain consequently using multiplication table and the equality $e^2=e$  that $p_8(e)=0$, $p_6(e)=p_7(e)=0$, $p_3(e)=0$ and $p_4(e)=0$, i.e. $e=0$.

Also, if $p_2(e)=0$, we have the equalities
$$
\begin{aligned}
p_3(e)&=p_3(e^2)=2(p_6(e)+p_7(e)-p_3(e))p_1(e)-(p_8(e)+p_5(e)-p_1(e))p_3(e),\\
p_4(e)&=p_4(e^2)=3(p_6(e)+p_7(e)-p_3(e))p_3(e)-(2p_8(e)-3p_1(e))p_4(e),\\
p_6(e)&=p_6(e^2)=(p_3(e)-p_6(e)-p_7(e))p_5(e)-(p_8(e)+p_5(e)-p_1(e))p_6(e),\\
p_7(e)&=p_7(e^2)=(p_6(e)+p_7(e)-p_3(e))p_8(e)-(p_8(e)+p_5(e)-p_1(e))p_7(e),
\end{aligned}
$$
which can be rewritten in the form
\begin{equation}\label{eqs}
\begin{aligned}
(p_8(e)+p_5(e)-p_1(e)+1)p_3(e)&=2(p_6(e)+p_7(e)-p_3(e))p_1(e),\\
(2p_8(e)-3p_1(e)+1)p_4(e)&=3(p_6(e)+p_7(e)-p_3(e))p_3(e),\\
(p_8(e)+p_5(e)-p_1(e)+1)p_6(e)&=(p_3(e)-p_6(e)-p_7(e))p_5(e),\\
(p_8(e)+p_5(e)-p_1(e)+1)p_7(e)&=(p_6(e)+p_7(e)-p_3(e))p_8(e).
\end{aligned}
\end{equation}

{\bf 2.} Assume that $p_2(e)=0$, $p_8(e)-p_5(e)-2p_1(e)\not=0$. Let us take
$a=\frac{p_3(e)-p_6(e)-p_7(e)}{p_8(e)-p_5(e)-2p_1(e)}$.
It is easy to see that $(p_3-p_6-p_7)\TT_a(e)=0$. So we may suppose that $(p_3-p_6-p_7)(e)=0$ initially.
Since $$2p_8(e)-3p_1(e)+1=2(p_8(e)+p_5(e)-p_1(e)+1)=2(p_8(e)-p_5(e)-2p_1(e))\not=0,$$ it follows from \eqref{eqs} that $p_3(e)=p_4(e)=p_6(e)=p_7(e)=0$. Let us denote $p_5(e)$ by $c$ and $p_8(e)$ by $d$. Then $p_1(e)=-(2c+1)$ and we have
$e=w_2(c,d)$.

{\bf 3.} Assume that $p_2(e)=0$, $p_8(e)-p_5(e)-2p_1(e)=0$. If $p_6(e)+p_7(e)-p_3(e)\not=0$, then it follows from \eqref{eqs} that $p_1(e)=p_5(e)=0$. This contradicts the equality $p_1(e)+2p_5(e)=-1$ proved above. So we have $p_3(e)=p_6(e)+p_7(e)$. Let us take $a=\frac{p_3(e)-4p_6(e)}{2}$. It is easy to see that $(4p_6-p_3)\TT_a(e)=0$. So we may suppose that $(4p_6-p_3)(e)=0$ initially. Let now consider two cases.

{\bf 3a.} Let $p_6(e)\not=0$. Then we have $p_6\SSS_{\frac{1}{p_6(e)}}(e)=1$ and so we may suppose that $p_6(e)=1$ initially. Then we have $p_3(e)=4$ and $p_7(e)=p_3(e)-p_6(e)=3$. Let us denote $p_5(e)$ by $-c$ and $p_4(e)$ by $d$. We have $p_1(e)=2c-1$ and $p_8(e)=p_5(e)+2p_1(e)=3c-2$. So
$e=w_3(c,d)$.

{\bf 3b.} Let $p_6(e)=0$. Then we have $p_3(e)=p_7(e)=0$. If $p_4(e)=0$, then we have
$e=w_2(c,d)$ for $c=p_5(e)$ and $d=p_8(e)$. Let now $p_4(e)\not=0$. There are $b\in\FF^*$ and $q\in\bar\FF$ such that $b^2p_4(e)=q$. If we denote $p_5(e)$ by $-c$, we obtain
$\SSS_b(e)=w_4(c,q)$.

Let us denote by $E$ the set $$\{w_1(c)\}_{c\in\FF}\cup\{w_2(c,d)\}_{c,d\in\FF}\cup \{w_3(c,d)\}_{c,d\in\FF}\cup\{w_4(c,q)\}_{c\in\FF, q\in\bar\FF}.$$
We prove that for any $w\in E$ and $f\in\Aut\big(W(2)\big)$ the condition $f(w)\in E$ is equivalent to $f(w)=w$ (it is clear that the union of sets listed in the theorem is disjoint in this case).

{\bf 1.} Suppose that $e=\SSS_b\TT_a(w_1(c))\in E$ ($a,c\in\FF$, $b\in\FF^*$). Since $p_2(e)=\frac{1}{b}$ and $p_2(w)\in\{0,1\}$ for any $w\in E$, we have $b=1$. Then we have $(p_1+2p_5)(e)=a-1$. Since $e\in E$ and $p_2(e)=1$, we have $(p_1+2p_5)(e)=-1$ and so $a=0$. We obtain that $e=w_1(c)$.

{\bf 2.} Suppose that $e=\SSS_b\TT_a(w_2(c,d))\in E$ ($a,c,d\in\FF$, $b\in\FF^*$). We have $p_3(e)=-4abc-2ab$, $p_6(e)=-abc$. We have $p_3(w)=4p_6(w)$ for any $w\in E$. Since $e\in E$, we have $-4abc-2ab=-4abc$, i.e. $a=0$. Then $e=\SSS_b(w_2(c,d))=w_2(c,d)$.

{\bf 3.} Suppose that $e=\SSS_b\TT_a(w_3(c,d))\in E$ ($a,c,d\in\FF$, $b\in\FF^*$). Analogously to the previous point we have $a=0$. Then $p_6(e)=b$. Since $p_6(w)\in\{0,1\}$ for any $w\in E$, we have $b=1$. Then $e=w_3(c,d)$.

{\bf 4.} Suppose that $e=\SSS_b\TT_a(w_4(c,q))\in E$ ($a,c\in\FF$, $b\in\FF^*$, $q\in\bar\FF$). Analogously to the second point of the present proof we have $a=0$. Then $e=-e_1-2e_8+c(3e_8-e_5+2e_1)+b^2qe_4$.
So $e\in E$ iff $b^2q=0$ or $b^2q\in\bar\FF$. Since $b\not=0$ and $q\not=0$, $b^2q\in\bar\FF$. But $q$ and $b^2q$ lie in the same coset of the group $(F^*)^2$ in the group $F^*$ and so by definition only one of them lies in $\bar\FF$ if $q\not=b^2q$.
Since $q$ lies in $\bar\FF$, we have $b^2q=q$ and so $e=-e_1-2e_8+c(3e_8-e_5+2e_1)+qe_4=w_4(c,q)$.
The theorem is proved.

\medskip

Let us introduce the algebras $W_1(t)$, $W_2(t,s)$, $W_3(t,s)$ and $W_4(t,u)$ for $s,t\in\FF$, $u\in\bar\FF$ by the following multiplication tables

\begin{center}
\begin{tabular}{llcll}
$W_1(t):$&
\begin{tabular}{c|c|c|}
&$v_1$&$v_2$\\ \hline
$v_1$&$-v_1+v_2$&$tv_2$\\ \hline
$v_2$&$0$&$0$\\ \hline
\end{tabular}&\phantom{xxxx}&
$W_2(t,s)$:&
\begin{tabular}{c|c|c|}
&$v_1$&$v_2$\\ \hline
$v_1$&$-v_1$&$tv_2$\\ \hline
$v_2$&$sv_2$&$0$\\ \hline
\end{tabular}\\
\\
$W_3(t,s)$:&
\begin{tabular}{c|c|c|}
&$v_1$&$v_2$\\ \hline
$v_1$&$-v_1$&$-v_2$\\ \hline
$v_2$&$tv_2-v_1$&$v_2+sv_1$\\ \hline
\end{tabular}&\phantom{xxxx}&
$W_4(t,u)$:&
\begin{tabular}{c|c|c|}
&$v_1$&$v_2$\\ \hline
$v_1$&$-v_1$&$-v_2$\\ \hline
$v_2$&$tv_2$&$uv_1$\\ \hline
\end{tabular}
\end{tabular}
\end{center}

It is easy to see that the algebra $W_1(t)$ corresponds to the element $w_1\left(\frac{t-2}{6}\right)\in W(2)$, the algebra $W_2(t,s)$ corresponds to the element $w_2\left(\frac{t+s-2}{6},\frac{t-s}{2}\right)\in W(2)$, the algebra $W_3(t,s)$ corresponds to the element $\SSS_{\frac{1}{6}}\left(w_3\left(\frac{3-t}{6},36s\right)\right)\in W(2)$, and the algebra $W_4(t,q)$ corresponds to the element $w_4\left(\frac{3-t}{6},q\right)\in W(2)$ (see the proof of Theorem \ref{idemp} for notation). We denote by $\mathcal{L}$ the set
$$\{W_1(t)\}_{t\in\FF}\cup\{W_2(t,s)\}_{t,s\in\FF}\cup \{W_3(t,s)\}_{t,s\in\FF}\cup\{W_4(t,u)\}_{t\in\FF, u\in\bar\FF}.$$

In the next corollary we show how to obtain a classification of two-dimensional algebras with a left quasi-unit using Theorem~\ref{idemp}.
Note that the classification of all two-dimensional algebras is known (see, for example,~\cite{GR11}).

\begin{corollary}\label{cor_quasiunit}
Let $R$ be a nonzero two-dimensional algebra with a left quasi-unit $e$. Then there is a unique algebra $W\in\mathcal{L}$ such that there is an isomorphism $\phi:R\rightarrow W$ satisfying the condition $\phi(e)=-v_1$.
\end{corollary}

{\bf Proof.} Let $\phi_0:R\rightarrow V_2$ be some isomorphism of vector spaces such that $\phi_0(e)=-v_1$. There is a unique multiplication $A$ on $V_2$ such that $\phi_0$ is an isomorphism from $R$ to $(V_2,A)$. Then $-v_1$ is a left quasi-unit of $(V_2,A)$, i.e. $A(v_1,A(x,y))-A(A(v_1,x),y)-A(x,A(v_1,y))=A(x,y)$ for all $x,y\in V_2$, and so $(V_2,A)$ is an idempotent of the algebra $W(2)$. By Theorem \ref{idemp} and the definition of the set $\mathcal{L}$ there are such $W\in\mathcal{L}$ and automorphism $f$ of $W(2)$ that $f$ sends the element of $W(2)$ corresponding to $(V_2,A)$ to the element of $W(2)$ corresponding to $W$. Then by Corollary \ref{isoaut} there is an isomorphism $\phi_1:(V_2,A)\rightarrow W$ such that $\phi_1(v_1)=v_1$. The map $\phi=\phi_1\phi_0$ satisfies all the required conditions. If we have another algebra $W'\in\mathcal{L}$ and an isomorphism $\phi':R\rightarrow W'$ such that $\phi'(e)=-v_1$, then there is an isomorphism $\phi'\phi^{-1}:W\rightarrow W'$ which sends $v_1$ to $v_1$. By Corollary \ref{isoaut} the elements of $W(2)$ corresponding to $W$ and $W'$ lie in the same orbit under the action of $\Aut\big(W(2)\big)$. But it follows from definition of $\mathcal{L}$ and Theorem \ref{idemp} that $W=W'$ in this case.

\appendix
\section*{Appendix A. Differences from the published version}
\addcontentsline{toc}{section}{Appendix A. Differences from the published version}

\noindent\textbf{Theorem~\ref{idemp} and Corollary~\ref{cor_quasiunit}:}   the proof has been corrected.

\end{document}